\documentclass[11pt]{article}

\usepackage{amsfonts,layout}

\newcommand{\R}{I\hspace{-1.5 mm}R}

\def\QED{\begin{flushright}
QED  \end{flushright}}

\newtheorem{Theorem}{Theorem}[section]

\newtheorem{Lemma}[Theorem]{Lemma}

\begin{document}

\title{A note on the regularity of solutions of Hamilton-Jacobi equations with superlinear growth in the gradient variable}
\author{P. Cardaliaguet\footnote{Laboratoire de Math\'ematiques,  UMR 6205, Universit\'e de Bretagne
Occidentale,  6 Av. Le Gorgeu,
BP 809, 29285 Brest. e-mail : Pierre.Cardaliaguet@univ-brest.fr. This work was partially supported by the ANR project MICA.}}

\maketitle

\section{Introduction}

We investigate the regularity of solutions to the Hamilton-Jacobi equation
\begin{equation}\label{HJ}
\left\{\begin{array}{ll}
-u_t(x,t)+b(x,t)|Du(x,t)|^q+f(x,t).Du(x,t)=0 & {\rm in}\; \R^N\times(0,T)\\
u(x,T)=g(x) & {\rm for}\; x\in \R^N
\end{array}\right.
\end{equation}
under the following assumptions: 
\begin{equation}\label{Hyp0}
q>1\;, 
\end{equation}
\begin{equation}\label{Hyp1}
\begin{array}{c}
\mbox{\rm  $b:\R^N\times(0,T)\to \R$, $f:\R^N\times(0,T)\to \R^N$ and $g:\R^N\to \R$}\\
\mbox{\rm   are continuous and bounded by some constant $M$},
\end{array}
\end{equation}
\begin{equation}\label{Hyp2}
b(x,t)\geq \delta >0 \qquad \forall (x,t)\in \R^N\times(0,T)
\end{equation}
for some $\delta>0$. 

Regularity of solutions of Hamilton-Jacobi equations with superlinear growth have been the object of several works (see 
in particular Lions \cite{Li}, Barles \cite{Ba1}, Rampazzo, Sartori \cite{RaSa}). 
Our aim is to show that $u$ is locally H\"{o}lder continuous with H\"{o}lder exponent and constant depending only $M$,
$\delta$, $q$ and $T$. What is new compared to the previous works is that the regularity does not depend
on the smoothness of the maps $b$, $f$ and $g$, but only on the growth condition. 
The motivation for this is the homogenization of Hamilton-Jacobi equations, where
such estimates are needed.  Here is our result.

\begin{Theorem}\label{main} There is some constant $\theta=\theta(M,\delta, q, T)$ and, for any $\tau>0$, some constant $K_\tau=K(\tau,M,\delta,q,T)$ 
such that, for any $x_0,x_1\in \R^N$, for any $t_0,t_1\in [0,T-\tau]$,  
$$
|u(x_0,t_0)-u(x_1,t_1)|\leq K_\tau \left( |x_0-x_1|^{(\theta-p)/(\theta-1)}+|t_0-t_1|^{(\theta-p)/\theta}\right)
$$
\end{Theorem}

The proof of the result relies on the representation of the solution $u$ of (\ref{HJ}) as the value function of a
problem of calculus of variations (see \cite{BaCa}, \cite{Ba}): Namely, setting $p=\frac{q}{q-1}$, we have
\begin{equation}\label{ControlPb}
u(x,t)=\inf\;  \left( \int_t^T a(x(s),s)\left|f(x(s),s)+x'(s)\right|^pds+g(x(T)) \right)
\end{equation}
where the infimum is taken over the set of functions $x(\cdot)\in W^{1,p}([t,T], \R^N)$ such that $x(t)=x$ and where 
$$
a(x,t)=\left(\frac{1}{b(x,t)}\right)^{p-1}\left(p^{-1/(p-1)}-p^{-p/(p-1)}\right)^{p-1}\;.
$$
From now on we work on the control representation of the solution $u$.
To simplify the notations, we assume without loss of generality that $b$ is also bounded by $M$ and satisfies
$$
a(x,t)\geq \delta >0 \qquad \forall (x,t)\in \R^N\times(0,T)
$$

The paper is organized as follows. In the first section, we use a kind of reverse H\"{o}lder inequality to prove that the optimal
solutions of (\ref{ControlPb}) are in some sense slightly ``more integrable" than what we could expect. In the second step we show
that this integrability implies the desired H\"{o}lder regularity for the value function. In Appendix, we prove the reverse H\"{o}lder inequality. \\

{\bf Aknowledgement : } We wish to thank Guy Barles for useful discussions.

\section{Estimate of the optimal of the controlled system}

The key remark of this section is Lemma \ref{ReguMorrey} stating that optimal controls are ``more integrable" than what could be expected. 
This is proved through several steps and the use of a reverse H\"{o}lder inequality.

\begin{Lemma}\label{BoundSol}
There is a constant $K\geq 0$ depending only on $M,\delta,p,T$, such that, for any 
optimal solution $\bar x$ of (\ref{ControlPb}) starting from $x_0$ at time $t_0$, we have
\begin{equation}\label{bound}
\int_{t_0}^{T}|\bar x'(s)|^pds \leq K\;.
\end{equation}
 \end{Lemma}

{\bf Proof of Lemma \ref{BoundSol} :  } Comparing $\bar x$ with the constant solution $\tilde x(t)=x_0$ we get
$$
\int_{t_0}^T a(\bar x(s),s)|f(\bar x (s),s)+\bar x'(s)|^pds+g(\bar x(T)) \leq \int_{t_0}^T a(x_0,s)|f(x_0,s)|^pds+g(x_0)
$$
with 
$$
g(x_0)-g(\bar x(T))\leq 2M\;,
$$
$$
\begin{array}{l}
\int_{t_0}^{T}a(\bar x(s),s)|f(\bar x (s),s)+\bar x'(s)|^pds \\
\begin{array}{rl}
\qquad  \geq &  \delta \int_{t_0}^{T}|f(\bar x (s),s)+\bar x'(s)|^pds\\
\geq &  \frac{\delta}{2^{p-1}} \left( \int_{t_0}^{T}|\bar x'(s)|^pds - M^p(T-t_0)\right) 
\end{array}
\end{array}
$$
and
$$
\int_{t_0}^{T}a(x_0,s)|f(x_0,s)|^pds \; \leq \; M^{p+1}(T-t_0)\;.
$$
Whence the result with $K=2^{p-1}(M^{p+1}T+2M)/\delta+M^pT$.
\QED

\begin{Lemma}\label{EstiOptiSol}
There are some constants $A\geq 1$ and $B\geq 0$ depending only on $M,\delta,p,T$, such that, for any 
optimal solution $\bar x$ of (\ref{ControlPb}) starting from $x_0$ at time $t_0$, we have
\begin{equation}\label{Ineqab}
\frac{1}{h}\int_{t_0}^{t_0+h}|\bar x'(s)|^pds \leq A\left(\frac{1}{h}\int_{t_0}^{t_0+h}|\bar x'(s)|ds\right)^p+B\qquad \forall h\in [0, T-t_0]
\end{equation}
 \end{Lemma}

{\bf Proof of Lemma \ref{EstiOptiSol} : } Let us fix $h\in (0, T-t_0)$ and let us define
$$
\tilde x (t)=\left\{\begin{array}{ll}
\frac{\bar x(t_0+h)-x_0}{h}\,(t-t_0)+x_0 & {\rm if }\; t\in [t_0,t_0+h]\\
\bar x (t) & {\rm otherwise}
\end{array}\right.
$$
Since $\bar x$ is optimal and $\tilde x(T)=\bar x(T)$ we have
$$
\int_{t_0}^{t_0+h}a(\bar x(s),s)|f(\bar x (s),s)+\bar x'(s)|^pds\leq 
\int_{t_0}^{t_0+h}a(\tilde x(s),s)|f(\tilde x (s),s)+\tilde x'(s)|^pds
$$
Then we get the desired result by noticing that
$$
\begin{array}{l}
\int_{t_0}^{t_0+h}a(\bar x(s),s)|f(\bar x (s),s)+\bar x'(s)|^pds \\
\begin{array}{rl}
\qquad  \geq &  \delta \int_{t_0}^{t_0+h}|f(\bar x (s),s)+\bar x'(s)|^pds\\
\geq &  \frac{\delta}{2^{p-1}} \left( \int_{t_0}^{t_0+h}|\bar x'(s)|^pds - M^ph\right) 
\end{array}
\end{array}
$$
and
$$
\begin{array}{l}
\int_{t_0}^{t_0+h}a(\tilde x(s),s)|f(\tilde x (s),s)+\tilde x'(s)|^pds\\
\begin{array}{rl}
\qquad \leq & M \int_{t_0}^{t_0+h}|f(\tilde x (s),s)+\tilde x'(s)|^pds \\
\leq & 2^{p-1} M(M^ph+ \int_{t_0}^{t_0+h}|(\bar x(h)-x_0)/h|^pds\\
\leq & 2^{p-1} M(M^ph+ h^{1-p} (\int_{t_0}^{t_0+h} |\bar x'(s)|ds)^p)
\end{array}
\end{array}
$$
\QED

In the following Lemma we get rid of the constant $B$ in (\ref{Ineqab}). Assume that $\alpha\in L^p([t_0,T],\R^+)$ satisfies
$$
\frac{1}{h}\int_{t_0}^{t_0+h}|\alpha(s)|^pds \leq A\left(\frac{1}{h}\int_{t_0}^{t_0+h}|\alpha(s)|ds\right)^p+B\qquad \forall h\in [0, T-t_0]
$$
Let $\displaystyle{  z(t)=\int_{t_0}^{t}\alpha(s)ds  }$ and 
$$
z_1(s)=\max\left\{z(s), (\frac{B}{A})^{1/p} (s-t_0)\right\}\qquad \forall s\in [t_0,T]\;.
$$
Set $\alpha_1(t)=z_1'(t)$. We note for later use that $z_1(t)\geq z(t)$ on $[t_0,T]$ and that, if $z_1(t)=z(t)$, then 
$\int_{t_0}^t (\alpha_1(s))^pds\leq \int_{t_0}^t (\alpha(s))^pds$. 
We claim:

\begin{Lemma}\label{2a}
$$
\frac{1}{h}\int_{t_0}^{t_0+h}|\alpha_1(s)|^pds \leq 2A\left(\frac{1}{h}\int_{t_0}^{t_0+h}|\alpha_1(s)|ds\right)^p\qquad \forall h\in [0, T-t_0]
$$
\end{Lemma}

{\bf Proof of Lemma \ref{2a} : } Let $\gamma= (B/A)^{1/p}$. If $z_1(t_0+h)=z(t_0+h)$, then from the definition of $z_1$
we have
$$
B\leq A (z(t_0+h)/h)^p
= A (z_1(t_0+h)/h)^p
$$
and therefore
$$
\int_{t_0}^{t_0+h} |\alpha_1(s)|^pds \leq  \int_{t_0}^{t_0+h} |\alpha(s)|^pds \; 
\leq \;  \frac{A}{h^{p-1}}(z(t_0+h))^p+Bh \;
\leq \; \frac{2A}{h^{p-1}}(z_1(t_0+h))^p
$$
If on the contrary $z_1(t_0+h)>z(t_0+h)$, then there is some $h_1<h$ such that $z_1(t_0+h_1)=z(t_0+h_1)$
and $z_1(s)=\gamma (s-t_0)$ on $[t_0+h_1,t_0+h]$. Then we have from the previous step
$$
\begin{array}{rl}
\int_{t_0}^{t_0+h} |\alpha_1(s)|^pds = & \int_{t_0}^{t_0+h_1} |\alpha_1(s)|^pds+\int_{t_0+h_1}^{t_0+h} |\alpha_1(s)|^pds\\
\leq & \frac{2A}{h_1^{p-1}} (z_1(t_0+h_1))^p+(h-h_1)\gamma^p\\
\leq  & 2A\gamma^p h_1+(h-h_1)\gamma^p\\
\leq & \frac{2A}{h^{p-1}}(z_1(t_0+h))^p
\end{array}
$$
 \QED

Next we show---in a kind of reverse H\"{o}lder inequality---that if a map satisfies the inequality given by Lemma \ref{2a}, then it is ``more integrable" than what we could expect.
There are several results of this nature in the literature since Gehring seminal work  \cite{Ge} (see for instance \cite{BeFr} and the references therein). 

\begin{Lemma}\label{RevHol}
Let $A>1$ and $p>1$. Then there are constants $\theta=\theta(A,p)>p$ and $C=C(A,p)>0$ such that, for any $\alpha\in L^p(0,1)$ such that
\begin{equation}\label{Hneq}
\frac{1}{h}\int_0^h |\alpha(s)|^pds \leq A \left(\frac{1}{h} \int_0^h |\alpha(s)|ds\right)^p \qquad \forall h\in [0,1],
\end{equation}
we have 
$$
\int_0^h |\alpha(s)|ds \leq C\|\alpha\|_p h^{1-1/\theta}\qquad \forall h\in [0,1]\;.
$$
Moreover, the optimal choice of  $ \theta$ is such that $\gamma =1-1/\theta$ is the smallest root of $\varphi(s) =s^p-A(1-p+ ps)$. 
\end{Lemma}

A possible proof of the Lemma is the following: using Gehring's result we can show that a map $\alpha$ satisfying (\ref{Hneq}) belongs in some $L^{r}$ 
for some $r>p$ with a $L^r$ norm controlled by its $L^p$ norm, and then use H\"older inequality. We have choosen to present in Appendix a new and direct proof 
using a completely different approach. \\

Combining Lemma \ref{EstiOptiSol}, Lemma \ref{2a} and Lemma \ref{RevHol} we get:

\begin{Lemma}\label{ReguMorrey}
There are constants $\theta>p$ and $C$ depending only on $M,\delta, p, T$ such that, for any $x_0\in \R^N$ and any $t_0<T$, if $\bar x$ is optimal for the initial position
$x_0$ at time $t_0$, then 
$$
\int_{t_0}^{t_0+h} |\bar x'(s)|ds\leq C(T-t_0)^{1/\theta-1/p} h^{1-1/\theta} \qquad \forall h\in [t_0, T]
$$
\end{Lemma}

{\bf Proof of Lemma \ref{ReguMorrey} : } Let $\bar x$ be optimal for $(x_0,t_0)$. From Lemma \ref{EstiOptiSol} we know that 
$$
\frac{1}{h}\int_{t_0}^{t_0+h}|\bar x'(s)|^pds \leq A\left(\frac{1}{h}\int_{t_0}^{t_0+h}|\bar x'(s)|ds\right)^p+B\qquad \forall h\in [0, T-t_0]
$$
for some constants $A,B$ depending only on $M, \delta,T$ and $p$. Setting $\alpha(t)=|\bar x'(t)|$, $z(t)=\int_{t_0}^{t}\alpha(s)ds$,
$z_1(t)=\max\{z(t), (B/A)^{1/p}\,(t-t_0)\}$ and $\alpha_1(t)=z_1'(t)$, we have from Lemma \ref{2a}:
$$
\frac{1}{h}\int_{t_0}^{t_0+h}(\alpha_1(s))^pds \leq 2A\left(\frac{1}{h}\int_{t_0}^{t_0+h}\alpha_1(s)ds\right)^p\qquad \forall h\in [0, T-t_0]
$$
Applying Lemma \ref{RevHol} to the constants $p$ and $2A$ and with a proper scalling, we get that there exists
$\theta>p$ and $C_2$ depending only on $M, \delta,T$ and $p$ such that
$$
\int_{t_0}^{t_0+h}|\bar x'(s)|ds\leq \int_{t_0}^{t_0+h}\alpha_1(s)ds \leq   (T-t_0)^{\frac{1}{\theta}-\frac{1}{p}}C_2\|\alpha_1\|_{p}h^{1-1/\theta}
$$
where we can estimate $\|\alpha_1\|_{L^p([t_0,T])}$ as follows: Let 
$$
\bar t=\max\{t\in [t_0,T]\;|\; z_1(t)=z(t)\}\;.
$$
Then
$$
\int_{t_0}^T \alpha_1^p(s)ds=  \int_{t_0}^{\bar t}\alpha_1^p(s)ds + \int_{\bar t}^T \left(\frac{B}{A}\right)^p ds\leq \int_{t_0}^{\bar t}\alpha^p(s)ds +\left(\frac{B}{A}\right)^p(T-\bar t)
$$
where, from Lemma  \ref{BoundSol}, we have
$$
\int_{t_0}^{\bar t}\alpha^p(s)ds\leq K\;.
$$
Therefore $\|\alpha_1\|_p\leq C_3$, where $C_3=C_3(M, \delta, p,T)$ and the proof is complete.
\QED

\section{Regularity of the value function}

We are now ready to prove Theorem \ref{main}.\\

{\bf Space regularity : } Let $x_0, x_1\in \R^N$, $t_0<T$. We assume that
\begin{equation}\label{hyph}
|x_1-x_0|\leq C\frac{(1-p/\theta)}{p-1}(T-t_0)^{1-1/p}\; \wedge \; 1\;,
\end{equation}
where $C$ and $\theta$ are the constants which appear in Lemma \ref{ReguMorrey}.
We claim that 
\begin{equation}\label{SpaceRegu}
u(x_1,t_0)-u(x_0,t_0)\leq  \; K_1 (T-t_0)^{-(p-1)(\theta-p)/(p(\theta-1))} |x_1-x_0|^{(\theta-p)/(\theta-1)}
\end{equation}
where $K_1=K_1(M, p, T, \delta)$. 

Indeed, let $\bar x$ be an optimal trajectory for $(x_0,t_0)$. For $h\in (0, T-t_0)$ let 
$$
\tilde x (t)=\left\{\begin{array}{ll}
\frac{\bar x(t_0+h)-x_1}{h}\,(t-t_0)+x_1 & {\rm if }\; t\in [t_0,t_0+h]\\
\bar x (t) & {\rm otherwise}
\end{array}\right.
$$
From Lemma \ref{ReguMorrey} we have
$$
|\bar x(t_0+h)-x_0|\leq \int_{t_0}^{t_0+h}|\bar x'(s)|ds\leq C(T-t_0)^{1/\theta-1/p}h^{1-1/\theta}\;.
$$
Therefore, since $\tilde x(T)=\bar x(T)$, we have
$$
\begin{array}{l}
u(x_1,t_0)\; \leq \; \int_{t_0}^{T} a(\tilde{x}(s),s)|f(\tilde x(s),s)+\tilde x'(s)|^pds+g(\tilde x(T)) \\
\begin{array}{rl}
\qquad \leq & u(x_0,t_0)+\int_{t_0}^{t_0+h} a(\tilde{x}(s),s)|f(\tilde x(s),s)+\tilde x'(s)|^pds\\
\leq & u(x_0,t_0)+M 2^{p-1} ( M^p h+h^{1-p}|\bar x(t_0+h)-x_1|^p)\\
 \leq & u(x_0,t_0)+M 2^{p-1} ( M^p h+h^{1-p}(|\bar x(t_0+h)-x_0|+|x_0-x_1|)^p\\
\leq &  u(x_0,t_0)+M 2^{p-1} ( M^p h+h^{1-p}(C_0h^{1-1/\theta}+|x_0-x_1|)^p\\
\end{array}
\end{array}
$$
where we have set $C_0=C(T-t_0)^{1/\theta-1/p}$. 
Choosing 
$$
h=\left(\frac{1}{C_0}\frac{p-1}{1-p/\theta}|x_0-x_1|\right)^{\theta/(\theta-1)}
$$
we have $h\leq (T-t_0)$ (from assumption (\ref{hyph})) and therefore
$$
u(x_1,t_0)-u(x_0,t_0)\; \leq \; K_1' (T-t_0)^{-(\theta-p)(p-1)/(p(\theta-1))} |x_1-x_0|^{(\theta-p)/(\theta-1)}
$$
where $K_1'=K_1'(M, \delta,p,T)$. Whence (\ref{SpaceRegu}).
\vspace{3mm}

{\bf Time regularity : } Let $x_0$ be fixed and $t_0<t_1< T-\tau$. We assume that
\begin{equation}\label{Condt0t1}
t_1-t_0\leq K_3 \tau^{(2\theta p-p-\theta)/ (p(\theta-1))}
\end{equation}
for some constant $K_3=K_3(M, \delta,p,T)$ to be fixed later, where 
$\theta$ is given by Lemma (\ref{ReguMorrey}). We claim that
$$
|u(x_0,t_0)-u(x_0,t_1)|\leq K_2\tau^{-(\theta-p)/\theta}(t_1-t_0)^{(\theta-p)/\theta}
$$
for some constant $K_2=K_2(M,\delta, p,T)$.

Indeed, let $\bar x$ be optimal for $(x_0,t_1)$. Then setting
$$
\tilde x (t)=\left\{\begin{array}{ll}
x_0 & {\rm if }\; t\in [t_0,t_1]\\
\bar x(t) & {\rm otherwise}
\end{array}\right.
$$
we have
$$
\begin{array}{rl}
u(x_0,t_0) \; \leq & \int_{t_0}^{t_1} a(x_0,s)|f(x_0,s)|^pds+ u(x_0,t_1)\\
\leq & M^{p+1}|t_1-t_0|+u(x_0,t_1)
\end{array}
$$
which gives the desired inequality provided $K_2$ is sufficiently large.

To get a reverse inequality, let $\bar x$ be now optimal for $(x_0,t_0)$.
Using Lemma \ref{ReguMorrey} we have that 
\begin{equation}\label{Estixt1}
\begin{array}{rl}
|\bar x(t_1)-x_0|\; \leq  & \int_{t_0}^{t_1}|\bar x'(s)|ds\; \leq \;  C(T-t_0)^{1/\theta-1/p}(t_1-t_0)^{1-1/\theta}\\
\leq & C\frac{(1-p/\theta)}{p-1}(T-t_1)^{1-1/p}\; \wedge \; 1 
\end{array}
\end{equation}
from the choice of $t_1-t_0$ in (\ref{Condt0t1}) and $K_3$ sufficiently small. Note that we have
$$
u(\bar x(t_1),t_1)\; \leq \;  u(x_0,t_0)-\int_{t_0}^{t_1} a(\bar x (s),s) |f(\bar x(s),s)+\bar x'(s)|^pds \leq u(x_0,t_0)
$$
Hence, using the space regularity of $u$ (recall that (\ref{Estixt1}) holds) we get
$$
\begin{array}{rl}
u(x_0,t_1) \; \leq & u(x_0,t_1)-u(\bar x(t_1),t_1)+u(x_0,t_0) \\
\leq &  u(x_0,t_0)+ K_1 C(T-t_1)^{-(p-1)(\theta-p)/(p(\theta-1))} |\bar x(t_1)-x_0|^{(\theta-p)/(\theta-1)} \\
\leq & u(x_0,t_0)+ K_2\tau^{-(\theta-p)/\theta}(t_1-t_0)^{(\theta-p)/\theta}
\end{array}
$$
\QED


\section{Appendix : proof of Lemma \ref{RevHol}}

 We note later use that the map $\varphi(s)= s^p-A(1-p+ps)$ has two roots, the smallest one---denoted by $\gamma$---belonging to the interval
$(1-1/p, A^{1/(p-1)})$, the other one being larger than $A^{1/(p-1)}$.
Moreover, if $\varphi(s)\leq 0$, then $s\geq \gamma$. 

Let 
$$
{\cal E}=\left\{\alpha\in L^p(0,1)\;,\; \alpha\geq 0, \; \alpha \; \mbox{\rm satisfies (\ref{Hneq}) and $\|\alpha\|_p\leq 1$}\right\}
$$
We note that ${\cal E}$ is convex, closed and bounded in $L^p(0,1)$. Therefore the problem
$$
\xi(\tau)=\max\left\{ \int_0^\tau \alpha(s)ds \; , \; \alpha \in {\cal E}\right\} 
$$
has a unique maximum denoted $\bar \alpha_\tau$ for any $\tau\in (0,1]$ (uniqueness comes from the fact that inequality (\ref{Hneq}) is positively homogeneous, which 
entails that at the optimum inequality $\|\alpha\|_p\leq 1$ is an equality). 

In order to prove the Lemma, we only need to show that 
\begin{equation}\label{IneqReduc}
\xi(\tau)\leq C \tau^{\gamma}\qquad \forall \tau\in [0,1]
\end{equation}
for a suitable choice of $C$, because again inequality (\ref{Hneq}) is positively homogeneous in $\alpha$. \\

The proof of (\ref{IneqReduc}) is achieved in two steps. In the first one, we explain the structure of the optima.
Then we deduce from this that $\xi$ satisfies a differential equation, which gives the desired bound.\\

{\bf Structure of the optima : } We claim that there is some $\bar \tau>0$ such that for any $\tau\in (0,\bar \tau)$, 
$$
\bar \alpha_\tau(t)=\left\{ \begin{array}{ll}
a_\tau & {\rm on }\; [0,\tau)\\
b_\tau & {\rm on }\; [\tau, \tau_1)\\
A^{-1/p} \gamma t^{\gamma-1} & {\rm on }\; [\tau_1, 1]
\end{array}\right.
$$
where $0<b_\tau\leq  a_\tau$ and $\tau < \tau_1< 1$.

{\it Proof of the claim : } Let $\bar x_\tau(t)=\int_0^t \bar \alpha_\tau(s)ds$. To show that $\bar \alpha_\tau$ is constant on $[0,\tau)$, we introduce
the map $\alpha (s)=\frac{\bar x_\tau(\tau)}{\tau}$ on $[0, \tau)$, $\alpha= \bar \alpha_\tau$ otherwise.
Then $\alpha$ belongs to ${\cal E}$ and is also optimal. Hence $\alpha=\bar \alpha_\tau$, which shows that $\bar \alpha_\tau$ is constant on $[0,\tau)$. 

With similar arguments we can prove that, if there is a strict inequality in (\ref{Hneq}) for $\bar \alpha_\tau$ at some  $h\geq \tau$, then $\bar \alpha_\tau$
is locally constant in a neighbourhood of $h$ in $[\tau,1]$. In particular, since $\bar \alpha_\tau$ is constant on $[0,\tau)$,
inequality (\ref{Hneq}) is strict for $\bar \alpha_\tau$ at $\tau$, and there is a maximal interval $[\tau,\tau_1)$ on which $\bar \alpha_\tau$ is constant.
We set $a_\tau=\bar \alpha_\tau(0^+)$ and $b_\tau=\bar \alpha_\tau(\tau^+)$. 

In order to show that $a_\tau\geq b_\tau$, we prove that
\begin{equation}\label{tata}
\mbox{\rm  the map $t\to \bar x_\tau(t)/t$ is nonincreasing.}
\end{equation}
 Indeed, let $t>0$ be fixed and
$x(s)= \max\{\bar x_\tau(s), \frac{\bar x_\tau(t)}{t} s\}$ if $s\in [0,t]$ and $x=\bar x_\tau$ otherwise. Let us check that $x'$ is admissible  and optimal. 
Let $I\subset (0,t)$ be the open set $\{x> \bar x_\tau\}$. We can write $I$ as the (at most) enumerable union of disjoint intervals $(c_i,d_i)$. 
Since $x$ is affine on each interval $(c_i,d_i)$ with $x(c_i)=\bar x(c_i)$ and $x(d_i)=\bar x(d_i)$ we have
\begin{equation}\label{tata1}
\int_{c_i}^{d_i} |x'|^p\leq \int_{c_i}^{d_i} |\bar x_\tau'|^p\qquad \forall i\;.
\end{equation}
Since $x'=\bar x_\tau$ a.e. in $[0,1]\backslash I$, we get $\|x'\|_p\leq \|\bar x_\tau'\|_p=1$. Moreover, for any $h>0$ such that $x(h)=\bar x_\tau(h)$,
(\ref{tata1}) and the admissibility of $\bar x_\tau'$ also give
$$
\frac{1}{h} \int_0^h |x'|^p \leq \frac{1}{h} \int_0^h |\bar x_\tau '|^p \leq A \left(\frac{1}{h} \bar x_\tau(h)\right)^p= A \left(\frac{1}{h}  x(h)\right)^p\;.
$$
If $x(h)>\bar x_\tau(h)$, let $h_1= \max\{s\leq h \; |\; x(s)=\bar x_\tau(s)\}$. Then $x(s)=\frac{\bar x_\tau(t)}{t}s$ on $[h_1, h]$ and so
$$
\begin{array}{rl}
\frac{1}{h} \int_0^h |x'|^p \; = & \frac{1}{h} \int_0^{h_1} |x'|^p + \frac{1}{h} \int_{h_1}^h |x'|^p\\
\leq &  \frac{Ah_1}{h} \left(\frac{1}{h_1}  x(h_1)\right)^p+  \frac{1}{h}(h-h_1) \left(\frac{\bar x_\tau(t)}{t}\right)^p\\
\leq & \frac{Ah_1}{h} \left(\frac{\bar x_\tau(t)}{t} \right)^p+  A \frac{1}{h}(h-h_1) \left(\frac{\bar x_\tau(t)}{t}\right)^p\\
\leq & A \left(\frac{\bar x_\tau(t)}{t} \right)^p  = A \left(\frac{1}{h} x(h)\right)^p
\end{array}
$$
So $x'$ is admissible. Since $x(\tau)\geq \bar x_\tau(\tau)$, $x$ is also optimal. So $x=\bar x_\tau$ and (\ref{tata}) is proved. 

Note that (\ref{tata}) implies that $a_\tau\geq b_\tau$ and
\begin{equation}\label{InIn}
\frac{s \bar x_\tau'(s)}{\bar x_\tau(s)} \leq 1< A^{1/(p-1)}\qquad \mbox{\rm for a.e. } s\in [0,1]\;.
\end{equation}

Let us now assume that $\tau_1<1$. To prove that $\bar \alpha_\tau(s)=A^{-1/p} \gamma s^{\gamma-1}$ on $[\tau_1,1]$, we show that there is an equality in (\ref{Hneq})
for $\bar \alpha_\tau$ on $[\tau_1, 1]$. Indeed, otherwise, $\bar \alpha_\tau$ is constant on some maximal interval $(u,v)$ with $\tau_1\leq u<v\leq 1$. 
We note that equality holds in (\ref{Hneq}) at $u$ because $\bar \alpha_\tau$ is not locally constant at this point. Taking the derivative with respect to $h$
in (\ref{Hneq}) at $u$ we get
$$
(\bar \alpha_\tau(u^+))^p\leq -\frac{(p-1)A}{u^{p}}(\bar x_\tau(u))^p+\frac{pA}{u^{p-1}}(\bar x_\tau(u))^{p-1}\bar \alpha_\tau (u^+)\;,
$$
i.e.,
$$
(\frac{u\bar \alpha_\tau(u^+)}{\bar x_\tau(u)})^p- A(1-p +p \frac{u\bar \alpha_\tau(u^+)}{\bar x_\tau(u)})\leq 0\;.
$$
From the analysis of $\varphi$, this implies that  
$$
\bar \alpha_\tau(u^+)\geq \gamma \frac{\bar x_\tau(u)}{u}\;.
$$
Let us define
$$
x(s)=\bar x_\tau(s) \; {\rm on }\; [0, u],\; x(s)=\frac{\bar x_\tau(u)}{u^\gamma} s^{\gamma} \; {\rm on }\; [u,v], \; x(s)=\frac{x(v^-)}{\bar x_\tau(v)} \bar x_\tau(s) \; {\rm on }\; [v,1]
$$
and $\alpha=x'$. Since $\alpha(u^+)=\gamma \frac{\bar x_\tau(u)}{u}\leq \bar \alpha_\tau(u^+)$, one easily checks that $x\leq \bar x_\tau$ and 
$\alpha \leq \bar \alpha_\tau$ on $[0,1]$. Moreover, a straightforward verification shows that $x$ satisfies (\ref{Hneq}). Hence $x$ is also optimal, which is impossible.
So there is an equality in (\ref{Hneq}) for $\bar \alpha_\tau$ on $[\tau_1, 1]$. Taking the derivative in this equality shows that $\bar \alpha_\tau$ solves
$$
\bar \alpha_\tau^p(s)=-\frac{pA}{s^{p-1}}(\bar x_\tau(s))^p+\frac{A}{s^p}(\bar x_\tau(s))^{p-1}\bar \alpha_\tau (s)\qquad {\rm on}\; [\tau_1,1]\;.
$$
From (\ref{InIn}) and the analysis of $\varphi$, this implies that $\bar x_\tau'(s)=\gamma \frac{\bar x_\tau(s)}{s}$ on $[\tau_1,1]$. 
Hence $\bar x_\tau(s)=Cs^\gamma$ for some constant $C$. Since there is an
equality in (\ref{Hneq}) at $h=1$ and since $\|\bar \alpha_\tau\|_p=1$, $1=A (\bar x_\tau(1))^p$ and therefore $C=A^{-1/p}$. 

Finally we have to show that $\tau_1<1$ for any $\tau\in (0,\bar \tau)$. Indeed, assume otherwise that $\tau_1=1$ for arbitrary small $\tau$.
Since $x(t)=A^{-1/p} t^\gamma$ is admissible, we have $a_\tau\tau \geq A^{-1/p} \tau^\gamma$. Hence $a_\tau\to+\infty$ as $\tau\to 0^+$. 
Moreover the constraint $\|\bar \alpha_\tau\|_p=1$ implies that $b_\tau$ is bounded when $\tau \to 0^+$. Hence, for any $k$ large, we can
find $\tau>0$ such that $a_\tau>k b_\tau$. Writing inequality (\ref{Hneq}) at $h=k\tau$ then gives
$$
a_\tau^p\tau \leq \frac{A}{(1+k)^{p-1}\tau^{p-1}} (a_\tau \tau+b_\tau k\tau )^p\leq \frac{A}{(1+k)^{p-1}\tau^{p-1}} (2 a_\tau \tau)^p=
\frac{2^pA}{(1+k)^{p-1}} a_\tau^p \tau
$$
whence a contradiction since $k$ is arbitrarily large. \\

{\bf A differential equation for $\xi$ : } To complete the proof of (\ref{IneqReduc}), we are going to show that $\xi$ is locally Lipschitz continous
and satisfies
\begin{equation}\label{de}
(-\tau)\xi'(\tau)+\gamma \xi(\tau)=0 \qquad \mbox{\rm for a.e. } \tau \in (0,\bar \tau)\;.
\end{equation}
From this (\ref{IneqReduc}) follows easily for a suitable choice of $C$. \\

{\it Proof of (\ref{de}) : } Let us extend the optimal solutions by $A^{-1/p}\gamma s^{\gamma-1}$ on $[1,+\infty)$ for $\tau\in (0,\bar \tau)$. 
For $\lambda>0$, let 
$$
\alpha_{\lambda\tau}(s)=\bar \alpha_\tau(\lambda s)\qquad s\geq 0\;.
$$
Then $\alpha_{\lambda\tau}$ satisfies (\ref{Hneq}) and 
$$
\|\alpha_{\lambda\tau}\|_p= \lambda^{-1/p}\left(1+\int_1^\lambda \bar \alpha_\tau^p\right)^{1/p}
$$
Hence $\alpha_{\lambda\tau}/\|\alpha_{\lambda\tau}\|_p$ is admissible and
\begin{equation}\label{xixi}
\xi\left(\frac{\tau}{\lambda}\right)\geq \frac{\int_0^{\tau/\lambda}\alpha_{\lambda\tau}}{\|\alpha_{\lambda\tau}\|_p}= 
\frac{\lambda^{1/p-1}\xi(\tau)}{\left(1+\int_1^\lambda \bar \alpha_\tau^p\right)^{1/p}}\;.
\end{equation}
with an equality for $\lambda =1$. In particular, this shows that $\xi$ is locally Lipschitz continuous in $(0,1]$.
Moreover, at each point $\tau$ at which $\xi$ has a derivative, we have, by taking the derivative with respect to $\lambda$
at $\lambda =1$ in (\ref{xixi}):
$$
(-\tau)\xi'(\tau)= (1/p-1) \xi(\tau)-\frac{\xi(\tau)}{p}\bar\alpha_\tau^p(1)=\xi(\tau)(1/p-1-A\gamma^p/p)=-\gamma \xi(\tau)
$$
on $(0, \bar \tau)$. Whence (\ref{de}).
\QED


\end{document}